\documentclass[journal]{IEEEtran}

\usepackage{myStyleIEEE}
\makenomenclature

\usepackage[inline]{enumitem}

\usepackage{amsmath}
\usepackage{amssymb}
\usepackage{diagbox}
\usepackage{listings}
\algblockdefx{FORALLP}{ENDFAP}[1]
  {\textbf{for all }#1 \textbf{do in parallel}}
  {\textbf{end for}}

\usepackage{alphalph}
\usepackage{etoolbox}
\usepackage{booktabs}
\usepackage{multirow}
\usepackage{xcolor}

\usepackage[inline]{enumitem}

\usepackage{bm}
\usepackage{varwidth}
\usepackage{mathtools}

\definecolor{color_ConvD}{rgb}{0.8588,0.3333,0.3333}
\definecolor{color_bilevel}{rgb}{0.333,0.8588,0.6039}
\definecolor{color_StochD}{rgb}{0.3216,0.6863,0.8980}

\usepackage{pgfplotstable}

\tikzset{pics/.cd,
  Neuron/.style={
    code={\ifnum#1=0
    \fill (-0.5,-0.5) rectangle (0.5,0.5);
    \else
    \fi
    }}}
\newcommand*{\ReadOutElement}[4]{%
    \pgfplotstablegetelem{#2}{[index]#3}\of{#1}%
    \let#4\pgfplotsretval
}

\patchcmd{\subequations}{\alph{equation}}{\alphalph{\value{equation}}}{}{}

\newcommand\ddfrac[2]{{\displaystyle\frac{\displaystyle #1}{\displaystyle #2}}}

\IEEEoverridecommandlockouts

\begin{document}

\title{Data-driven Adaptive Benders Decomposition for the Stochastic Unit Commitment Problem}

\renewcommand{\theenumi}{\alph{enumi}}

\author{Baudouin~Vandenbussche, 
        Stefanos Delikaraoglou,
        Ignacio Blanco
        and~Gabriela~Hug
        }

\maketitle
\IEEEpeerreviewmaketitle

\begin{abstract}
This paper proposes a data-driven version of the Benders decomposition algorithm applied to the stochastic unit commitment (SUC) problem.
The proposed methodology aims at finding a trade-off between the size of the Benders master problem and the number of iterations until convergence.
Using clustering techniques, we exploit the information contained in the Lagrange multipliers of the Benders subproblems in order to aggregate the optimality cuts, without compromising the critical information that is passed to the master problem. In addition, we develop an outer parallelization scheme that finds the optimal solution of the  SUC problem by solving a series of less computationally intensive SUC instances for certain partitions of the scenario set.
Our computational results on the IEEE 3-Area RTS-96 power system, illustrate the improved performance of our data-driven Benders algorithm, in terms of solution time and problem size, compared both to the SUC extensive formulation and to the prevailing single- and multi-cut Benders formulations.

\end{abstract}



\vspace{-0cm}

\vspace{-8pt}

\section{Introduction} \label{sec:Intro}

The appeal for the reduction of greenhouse gas emissions and the establishment of a sustainable power system have prompted the widespread development of renewable energy sources. 
According to the annual report of REN21 \cite{REN21}, 57 countries have established 100\% renewable energy targets and the annual growth rate of renewables over the past decade is more than 5\%.
However, the transition towards a fully-renewable electricity sector requires a paradigm shift in the market design and the operation of the power system,
as the power production of wind turbines or solar PVs depends on weather conditions that can be only partially predicted ahead of real-time operation. 
On the contrary, conventional generators which serve the base load and provide the necessary flexibility to cope with the inherent variability and the uncertainty induced by the forecast errors of stochastic renewables, have to be committed well in advance of the actual operation, due to their technical requirements.

The need for improved decision-support tools that embrace a probabilistic view of system uncertainties has been drawing extensive attention in the recent literature. 
In this line, the works in \cite{tuohy2009unit, bouffard2005market, wang2008security} among others, propose a variant of the classical unit commitment (UC) problem based on stochastic optimization \cite{birge2011introduction}. The so-called stochastic unit commitment (SUC) co-optimizes day-ahead and real-time operations in order to minimize expected system cost, having a probabilistic description of uncertainty in the form of scenarios. 
To obtain though a reliable stochastic solution, we need to model accurately the predictive densities of renewable generation using a large number of scenarios. Considering that the traditional UC model is already a computationally challenging mixed-integer linear program (MILP), the SUC model with a large scenario set may result in an intractable optimization problem.

Two main approaches have been proposed to cope with this high computational burden. The first employs scenario reduction techniques \cite{heitsch, Bruninx} to approximate with reasonable accuracy the probability distribution of the underlying stochastic process using a smaller scenario set.
The second approach relies on decomposition schemes that exploit the special structure of the SUC problem. Lagrangian relaxation is applied to SUC in \cite{Papa_lagrangian}, while \cite{aravena2015distributed} proposes a distributed asynchronous version of this algorithm. The Progressive Hedging algorithm is employed in \cite{Ordoudis} and \cite{ryan2013toward} to solve the SUC problem, whereas the column-and-constraint generation method is used in \cite{Blanco}. The authors in \cite{Zheng2013} and \cite{wu2010accelerating} employ a Benders decomposition algorithm to improve the computational tractability of the SUC model when applied to large-scale power systems.

In this work, we combine the aforementioned approaches in order to develop a data-driven version of the Benders decomposition algorithm. 
Using clustering techniques, we exploit the information contained in the Lagrange multipliers of the second-stage subproblems, which reflect the sensitivities of the recourse actions with respect to first-stage decisions, in order to control efficiently the amount of data that is passed to the Benders master problem.
In addition, we propose an outer parallelization algorithm that leverages the statistical properties of the scenario set to solve the SUC problem in a two-step process and enables the use of high performance computing resources.
 In the first step, we partition the scenario set into smaller subsets and we solve in a distributed fashion the corresponding SUC problems. These SUC models are less computationally demanding than the original SUC due to the reduced number of scenarios involved. 
The full scenario set is then fed to a SUC instance in which the commitment variables that were common among all subsets are treated as fixed parameters to reduce the branching effort of the MILP solver.
The proposed data-driven Benders algorithm brings two major advantages: it improves the computational tractability of the SUC problem in terms of solution time and size of the problem, while it keeps intact the uncertainty information.

The rest of the paper is organized as follows. Section \ref{sec:mathmodel} provides the mathematical formulation of the extensive form of the SUC problem and its reformulation using Benders decomposition. Section \ref{sec:Improvements} presents the proposed data-driven Benders decomposition improvements and the outer parallelization scheme. Section \ref{sec:Results} discusses the parameter tuning for enhancing the performance of the proposed methods and analyses the results from the case study based on the standard IEEE RTS-96 power system. 
Finally, Section \ref{sec:conclusion} concludes the paper and provides suggestions for future work.


\section{Mathematical formulation} 
\label{sec:mathmodel}

In this section, we first provide the mathematical formulation of the stochastic unit commitment problem and we then explain how it can be decomposed according to the single- and multi-cut versions of the Benders decomposition algorithm.

\vspace{-8pt}

\subsection{Stochastic Unit Commitment}

The stochastic unit commitment model is formulated as the following two-stage stochastic programming problem:
\begin{subequations}
\label{mod:SUC}
\allowdisplaybreaks
\begin{flalign}
& \hspace{-0.4cm} \underset{\Xi^{\text{DA}},\Xi^{\text{BA}}}{\text{Minimize}} \;\; C^{\text{DA}} + C^{\text{BA}} =  \nonumber \\
&\hspace{-0.4cm} \sum_{t \in \mathcal{T}}\sum_{g \in \mathcal{G}}\left(C_{g}p_{gt}+C^{\mathrm{SU}}_{g}y_{gt}+C^{\mathrm{U}}_{g}r^{+}_{gt}+C^{\mathrm{D}}_{g}r^{-}_{gt}\right)+ \nonumber \\
& \hspace{-0.4cm} \sum_{\omega \in \Omega} \pi_{\omega}\sum_{t \in \mathcal{T}}\Big[\sum_{g \in \mathcal{G}}(C^{+}_{g}p^{+}_{gt\omega}-C^{-}_{g}p^{-}_{gt\omega})+C^{\text{shed}}\sum_{n \in \mathcal{N}}l^{\text{shed}}_{nt\omega}\Big]\label{eq:objfunct}
\end{flalign} 
subject to
\begin{align}
&u_{gt}=u^{0}_{g}, \quad \forall g,\forall t \leq L^{\mathrm{U}}_{g}+L^{\mathrm{D}}_{g},   \label{eq:initial_gen_state}\\ 
&\sum_{t'=t-UT_{g}+1}^{t}y_{g t'} \leq u_{gt}, \quad \forall g,\forall t > L^{\text{U}}_{g}+L^{\text{D}}_{g}, \label{eq:on_startup}\\
&\sum_{t'=t-DT_{g}+1}^{t}z_{gt'} \leq 1-u_{gt}, \quad \forall g,\forall t > L^{\mathrm{U}}_{g}+L^{\mathrm{D}}_{g}  ,
\label{eq:off_shtdwn}\\
&y_{gt}-z_{gt}=u_{gt}-u^{0}_{g}, \quad \forall g, \; t=0  ,
\label{eq:init_stup_shtdwn}\\
&y_{gt}-z_{gt}=u_{gt}-u_{gt-1}, \quad \forall g,\forall t \in \{1,..,\mathcal{T}\},  
\label{eq:chge_state}\\
&y_{gt}+z_{gt} \leq 1, \quad \forall g, t  ,
\label{eq:chge_state_not_sm_time}\\
&p_{gt}+r^{+}_{gt}\leq p_{gt-1}+r^{+}_{gt-1}+R^{\mathrm{U}}_{g}(u_{gt-1}+y_{gt}),  \; \forall g, t ,
\label{eq:rmp_up}\\
&p_{gt-1}-r^{-}_{gt-1}\leq p_{gt}-r^{-}_{gt}+R^{\mathrm{D}}_{g}(u_{gt}+z_{gt}),\quad \forall g, t,
\label{eq:rmp_dwn}\\
&p_{gt}-r^{-}_{gt}\geq  \underline{P}_g u_{gt}, \quad \forall g, t,
\label{eq:pmin}\\
&p_{gt}+r^{+}_{gt}\leq  \overline{P}_g (u_{gt}-z_{gt+1})+R^{\mathrm{D}}_{g}z_{gt+1}, \quad \forall g, t,
\label{eq:pmax}\\
&r^{+}_{gt}\leq R^{+}_{g},\quad \forall g, t,
\label{eq:reserve_up}\\
&r^{-}_{gt}\leq R^{-}_{g}, \quad \forall g, t,
\label{eq:reserve_dwn}\\
& w_{jt}\leq W_{j}, \quad \forall j, t,
\label{eq:wind_lim}\\
&\sum_{g \in \mathcal{G}_{n}}p_{gt}+ \!\!\! \sum_{j \in \mathcal{W}_{n}} \!\!\! w_{jt}= 
L_{nt}+ \!\!\! \sum_{\ell \in \mathcal{F}_{n}^{\text{in}} }\hat{f}_{\ell t}-\sum_{\ell \in \mathcal{F}_{n}^{\text{out}} }\hat{f}_{\ell t}, \forall n, t, \label{eq:power_balance} \\
&-F_{\ell} \leq \hat{f}_{\ell t}  \leq F_{\ell}, \quad \forall \ell, t,
\label{eq:flow_cap}\\
&\hat{f}_{\ell t}=B_{\ell}\sum_{n \in \mathcal{N}_{\ell}}\hat{\delta}_{nt}, \quad \forall \ell, t,
\label{eq:flow_val}\\
&\sum_{g \in \mathcal{G}_{n}}\left(p^{+}_{gt\omega}-p^{-}_{gt\omega}\right)\! + \!\! \sum_{j \in \mathcal{W}_{n}} \!\! \left(W^{*}_{jt\omega}-w_{jt}-w^{\mathrm{spill}}_{jtw}\right) +l^{\text{shed}}_{nt\omega}= \nonumber \\
&\!\! \sum_{\ell \in \mathcal{F}_{n}^{\text{in}} }\left(\tilde{f}_{\ell t\omega}-\hat{f}_{\ell t}\right)
- \!\!\! \sum_{\ell \in \mathcal{F}_{n}^{\text{out}} }\left(\tilde{f}_{\ell t\omega}-\hat{f}_{\ell t}\right), \quad \forall n, t, \omega,
\label{eq:2nd_stage_power_bal}\\
&p^{+}_{gt\omega} \leq r^{+}_{gt}, \quad \forall g, t, \omega,
\label{eq:reserve_lim_scenar_up}\\
&p^{-}_{gt\omega} \leq r^{-}_{gt}, \quad \forall g,t,\omega,
\label{eq:reserve_lim_scenar_dwn}\\
&w^{\mathrm{spill}}_{jt\omega} \leq W^{*}_{jt\omega}, \quad \forall j, t,\omega,
\label{eq:wspill}\\
&l^{\text{shed}}_{nt\omega} \leq L_{nt}, \quad \forall n, t, \omega,
\label{eq:lshed}\\
&-F_{\ell} \leq \tilde{f}_{\ell t\omega}  \leq F_{\ell}, \quad \forall \ell, t, \omega,
\label{eq:flow_cap_2nd_stage}\\
&\tilde{f}_{\ell t\omega}=B_{\ell}\sum_{n \in \mathcal{N}_{\ell}}\tilde{\delta}_{n t\omega}, \quad \forall \ell, t, \omega,
\label{eq:flow_val_2nd_stage} \\
& u_{gt} = \{0,1\}, \; y_{gt}= \{0,1\} \;,z_{gt}= \{0,1\}, \;\forall g,t,   \nonumber \\
& p_{gt},  r^{+}_{gt},  r^{-}_{gt} \ge 0, \forall g,t; \; \; w_{jt} \ge 0,\forall j,t  , \nonumber \\
&\hat{\delta}_{nt} \; \text{free}, \forall n,t; \; \hat{f}_{\ell t} \; \text{free}, \; \forall \ell,t,  \label{eq:vardecDA} \\
& p^{+}_{gt\omega},p^{-}_{gt\omega} \ge 0, \; \forall g,t,\omega,  \; w^{\mathrm{spill}}_{jt\omega} \ge 0, \; \forall j,t , \nonumber \\
&  \tilde{f}_{\ell t \omega} \; \text{free}, \; \forall \ell,t, \omega,  l^{\text{shed}}_{nt\omega}, \tilde{\delta}_{nt\omega} \; \text{free}, \forall n,t, \omega, \label{eq:vardecRT}
\end{align}
\end{subequations}
where $\Xi^{\text{DA}}$ = \{$u_{gt}, y_{gt},z_{gt},p_{gt}, r^{+}_{gt}, r^{-}_{gt},  \;\forall g,t; w_{jt}, \forall j,t;$ $\hat{\delta}_{nt}, \forall n,t ; \hat{f}_{\ell t}, \forall \ell, t $\} is the set of first-stage (day-ahead) decision variables and $\Xi^{\text{BA}}$ = \{$p^{+}_{gt\omega}$, $p^{-}_{gt\omega}, \forall g,t, \omega; w^{\text{spill}}_{jt\omega}, \forall j,t,\omega; \tilde{\delta}_{nt\omega}$, $l^{\text{shed}}_{nt\omega}, \forall n,t, \omega; \tilde{f}_{\ell t\omega}, \forall \ell,t, \omega $\} is the set of second-stage (balancing) decision variables.

The objective function \eqref{eq:objfunct} to be minimized is the total expected system cost that comprises both day-ahead ($C^{\text{DA}}$) and balancing cost ($C^{\text{BA}}$) components. The day-ahead part of \eqref{eq:objfunct} accounts for the energy production and start-up costs of all conventional units, denoted by $C_g$ and $C^{\text{SU}}_{g}$, respectively, as well as for the upward and downward reserve capacity procurement costs $C^{\text{U}}_g$ and $C^{\text{D}}_g$. The real-time cost component includes the re-dispatch cost  for every uncertainty realization $\omega \in \Omega$, namely: the upward and downward reserve deployment at the corresponding offer prices $C^{+}_{g}$ and $C^{-}_{g}$, as well as the involuntary load shedding at the penalty cost of $C^{\text{shed}}$.

The first-stage constraints \eqref{eq:initial_gen_state}-\eqref{eq:flow_val} ensure that the day-ahead schedule respects the technical limits of the power system. 
Using the binary variables $u_{gt}, y_{gt}, z_{gt}$, the initial state of the units in the beginning of the scheduling horizon as well as the resulting allowable start-up and shut-down actions are imposed through constraints \eqref{eq:initial_gen_state} and \eqref{eq:on_startup} - \eqref{eq:off_shtdwn}, respectively. 
Parameter $u^{0}_{g}$ denotes the initial commitment status of unit $g$, whereas parameter $L^{\text{U}}_g$ ($L^{\text{D}}_g$) indicates the number of time periods for which generator $g$ must be online (offline) from the beginning of the scheduling horizon. Parameters $UT_g$ and $DT_g$ denote the minimum up and down time for generator $g$.
Constraints \eqref{eq:init_stup_shtdwn} and \eqref{eq:chge_state} model the transition from start-up to shut-down state, while constraint \eqref{eq:chge_state_not_sm_time} states that unit $g$ can either start-up or shut-down at time period $t$.
 The upward and downward ramping limits $R^{\text{U}}_g$ and $R^{\text{D}}_g$ are enforced by constraints \eqref{eq:rmp_up} and \eqref{eq:rmp_dwn}, respectively, taking into account the energy production schedule $p_{gt}$ as well as the amount of upward and downward reserve capacity, denoted by $r^{+}_{gt}$ and $r^{-}_{gt}$, procured from each generator. Similarly, the minimum $\underline{P}_g$ and maximum $\overline{P}_g$ generation bounds are enforced by constraints \eqref{eq:pmin} and \eqref{eq:pmax}, where the maximum power to shut down the units is equal to $R^{\mathrm{D}}_{g}$. The procurement of upward and downward reserves is limited by the corresponding capacity offers $R^{+}_{g}$ and $R^{-}_{g}$ using constraints \eqref{eq:reserve_up} and \eqref{eq:reserve_dwn}. Constraint \eqref{eq:wind_lim} bounds wind power dispatch $w_{jt}$ to the installed wind power capacity $W_{j}$. Finally, the nodal power balance is enforced by the equality constraint \eqref{eq:power_balance} taking into account the day-ahead network power flows $\hat{f}_{\ell}$ which are restricted by the transmission capacity limits $F_{\ell}$ in constraints \eqref{eq:flow_cap} and \eqref{eq:flow_val} based on a DC flow approximation.

The second-stage constraints \eqref{eq:2nd_stage_power_bal} - \eqref{eq:flow_val_2nd_stage} model the balancing recourse actions for each wind power realization $\omega$. Constraints \eqref{eq:reserve_lim_scenar_up} and \eqref{eq:reserve_lim_scenar_dwn} ensure that the deployment of upward $p^{+}_{gt \omega}$ and downward $p^{-}_{gt \omega}$ reserves, respectively, does not exceed the corresponding procured quantities at the day-ahead stage. The amount of wind power production $w^{\text{spill}}_{jt\omega}$ that can be spilled as well as the allowable load shedding $l^{\text{shed}}_{nt\omega}$ at each node are bounded to the realized wind power production $W^{*}_{jt\omega}$ in each scenario $\omega$ and to the nodal demand $L_{nt}$ in constraints \eqref{eq:wspill} and \eqref{eq:lshed}, respectively. Equation \eqref{eq:2nd_stage_power_bal} ensures that conventional generation, wind power production and load are properly re-dispatched  such that the whole system zone remains in balance, whereas constraints \eqref{eq:flow_cap_2nd_stage} and \eqref{eq:flow_val_2nd_stage} imposes the transmission capacity limits on the real-time power flows $\tilde{f}_{\ell t \omega}$. Constraints \eqref{eq:vardecDA} and \eqref{eq:vardecRT}  are variable declarations. 


\vspace{-8pt}
\subsection{Benders decomposition algorithm}

The solution of the stochastic unit commitment model \eqref{mod:SUC} can become very computationally intensive when this mixed-integer, NP-hard problem is applied to large-scale power systems in combination with a large set of scenarios to describe accurately the wind power uncertainty. Nonetheless, exploiting the structure of the problem at hand, we can apply the Benders decomposition algorithm in order to reduce the involved computational burden, while ensuring also convergence to the global optimal solution. 
Below, we present the two main implementations of the Benders algorithm, namely, the single- and multi-cut versions, which we will use as foundation for the decomposition strategies proposed in this work.

To facilitate the exposition, we first present the multi-cut version of Benders adapted to two-stage stochastic optimization problems \cite{L-shape} and we then outline its main differences compared to the single-cut approach \cite{birge2011introduction}. In the multi-cut Benders algorithm, the extensive form \eqref{mod:SUC} of the SUC model is decomposed into a master problem that involves only first-stage  variables and into a set of subproblems, one per scenario $\omega$, that model the recourse actions. The Benders master problem in iteration $\nu$ is formulated as:
\begin{subequations}
\label{mod:Benders-Master-MC}
\begin{align}
& \underset{\Xi^{\text{DA}},\theta_{\omega}}{\text{Minimize}} \;\; \mathcal{Q}^{\text{MC}} =  C^{\text{DA}} + \sum_{\omega \in \Omega} \pi_{\omega} \theta_{\omega}
\end{align}
\vspace{-15pt}

subject to

\vspace{-8pt}

\begin{align}
& \text{constraints }  \eqref{eq:initial_gen_state}-\eqref{eq:flow_val}, \; \eqref{eq:vardecDA}, \label{eq:UC-FScns} \\
& \theta_{\omega} \ge \theta^{\text{min}}, \quad \forall \omega, \label{eq:MC-theta-lb}\\
& \theta_{\omega} \ge \Theta^{(k)}_{\omega}, \quad \forall k=\{1,...,\nu-1\}, \forall \omega,  \label{eq:MC-theta-lb2} \\
& \begin{aligned}
& \Theta^{(k)}_{\omega}  =  \mathcal{Q}^{(k)}_{\omega} + \sum_{t \in \mathcal{T}}\Big[   \sum_{g \in \mathcal{G}}\lambda^{+;(k)}_{gt\omega}(r^{+}_{gt}-r^{+;(k)}_{gt}) +\\
& \sum_{g \in \mathcal{G}}\lambda^{-;(k)}_{gt\omega}(r^{-}_{gt}-r^{-;(k)}_{gt})+\sum_{j \in \mathcal{W}}\lambda^{\text{W};(k)}_{jt\omega}(w_{jt}-w^{(k)}_{jt}) + \\
& \sum_{\ell \in \mathcal{L}}\lambda^{\text{F};(k)}_{\ell t\omega}(\hat{f}_{\ell t}-\hat{f}^{(k)}_{\ell t})\Big], \quad \forall k=\{1,...,\nu-1\}, \forall \omega. \label{eq:MC-cuts}
\end{aligned} 
\end{align}
\end{subequations}
The first-stage constraints of the SUC problem are included in \eqref{eq:UC-FScns}, whereas information about the second stage is conveyed to the master problem via the auxiliary variable $\theta_{\omega}$ that is bounded from below by  parameter $\theta^{\text{min}}$ in \eqref{eq:MC-theta-lb} and the set of optimality cuts \eqref{eq:MC-cuts}. 
These cuts are essentially supporting hyperplanes of a function that maps first-stage decisions to optimal recourse actions for each scenario $\omega$. At each iteration $\nu$ of the Benders algorithm, a new set of cuts \eqref{eq:MC-cuts} is included in the master problem, using the dual variables $\lambda$ of the second-stage subproblems that are formulated, according to \cite{conejo2006decomposition}, for every scenario $\omega' \in \Omega$ as:
\begin{subequations}
\label{mod:suprolem}
\begin{flalign}
& \underset{\Xi^{\text{S}}_{\omega'}}{\text{Minimize}} \quad \mathcal{Q}^{(\nu)}_{\omega'}= \nonumber \\ 
& \sum_{t \in \mathcal{T}}\Big[\sum_{g \in \mathcal{G}}(C^{+}_{g}p^{+}_{gt\omega'}-C^{-}_{g}p^{-}_{gt\omega'})
+ C^{\mathrm{sh}}\sum_{n \in \mathcal{N}}l^\mathrm{shed}_{nt\omega'}\Big] \label{eq:obj-MC} 
\end{flalign}
subject to
\begin{flalign}
 &  \text{constraints }  \eqref{eq:2nd_stage_power_bal}-\eqref{eq:flow_val_2nd_stage}, \; \eqref{eq:vardecRT}, \quad \omega=\omega' , \label{eq:UC-SScns} \\
& r^{+}_{gt}=r^{+; (\nu)}_{gt} : \lambda^{+; (\nu)}_{gt\omega'}, \quad \forall g,t, \label{eq:FixCns1}\\
& r^{-}_{gt}=r^{-; (\nu)}_{gt} : \lambda^{-; (\nu)}_{gt\omega'}, \quad \forall g,t,\\
& w_{jt}=w^{(\nu)}_{jt} \;\; : \lambda^{\text{W}; (\nu)}_{jt\omega'}, \quad \forall j,t,\\
& \hat{f}_{\ell t}=\hat{f}^{(\nu)}_{\ell t} \;\;\;  : \lambda^{\text{F}; (\nu)}_{\ell t\omega'}, \quad \forall \ell,t  , \label{eq:FixCns4}
\end{flalign}
\end{subequations}
where 
$\Xi^{\text{S}}_{\omega'}$ = $\Xi^{\text{BA}}_{\omega'}$ $\bigcup$ $ \{r^{+}_{gt}, r^{-}_{gt}, \forall g,t; w_{jt}, \forall j,t ; \hat{f}_{\ell t}, \forall \ell,t\}$ is the set of primal optimization variables of the Benders subproblem for scenario $\omega'$, and $\lambda^{+}_{gt\omega'},\;  \lambda^{-}_{gt\omega'},\; \lambda^{\text{W}}_{jt\omega'}, \; \lambda^{\text{F}}_{\ell t\omega'}$ are the Lagrange multipliers of the `fixing' constraints \eqref{eq:FixCns1} - \eqref{eq:FixCns4} associated with the optimal solution of master problem at iteration $\nu$. It should be noted that allowing for load shedding and wind spillage during real-time operation, the second-stage constraints \eqref{eq:UC-SScns} are always satisfied and thus it is not necessary to include feasibility cuts in the master problem \eqref{mod:Benders-Master-MC}. 
At the end of every iteration $\nu$, an upper and a lower bound, denoted as $\overline{\xi}^{(\nu)}$ and $ \underline{\xi}^{(\nu)}$, respectively, are calculated as:
\begin{align}
    \overline{\xi}^{(\nu)} = C^{\text{DA}} + \sum_{\omega}\pi_{\omega} \mathcal{Q}^{(\nu)}_{\omega}, \quad \underline{\xi}^{(\nu)} = \mathcal{Q}^{\text{MC}} \label{eq:LB}
\end{align}
\vspace{-8pt}

\noindent and the Benders algorithm terminates when $\overline{\xi}^{\nu} - \underline{\xi}^{\nu} \le \varepsilon$, where $\varepsilon$ is a user-defined convergence threshold.

The single-cut version of Benders decomposition as presented in \cite{L-shape} adds only one cut per iteration in the master problem, which in turn is formulated as:
\begin{subequations}
\label{mod:Benders-Master-SC}
\begin{align}
& \underset{\Xi^{\text{DA}},\theta}{\text{Minimize}} \;\; \mathcal{Q}^{\text{SC}} =  C^{\text{DA}} + \theta \label{eq:obj-SC}
\end{align}
subject to
\begin{flalign}
& \text{constraints }  \eqref{eq:initial_gen_state}-\eqref{eq:flow_val}, \; \eqref{eq:vardecDA}, \; \eqref{eq:MC-cuts},  \\
& \theta \ge \theta^{\text{min}}, \\
& \theta \ge \sum_{\omega \in \Omega} \pi_{\omega} \Theta^{(k)}_{\omega}, \quad \forall k=\{1,...,\nu-1\}  .
\end{flalign}
\end{subequations}

The convergence criterion is the same as in the multi-cut version, albeit the lower bound in \eqref{eq:LB} is now calculated  as $\underline{\xi}^{(\nu)} = \mathcal{Q}^{\text{SC}}$. Comparing the two variations of the Benders algorithm, one notices that the size of the master problem \eqref{mod:Benders-Master-MC} in the multi-cut version grows significantly faster compared to its single-cut counterpart \eqref{mod:Benders-Master-SC}. 
By adding $|\Omega|$ new constraints (cuts) in every iteration $\nu$, where $|\Omega|$ is the cardinality of the scenario set $\Omega$, the multi-cut approach preserves the complete information that is passed from the sub-problems to the master problem and enables convergence to the optimal solution in a lower number of iterations \cite{BIRGE}. However, this comes at the expense of larger master problem instances and higher memory requirements. On the other hand, the single-cut version aggregates all the second-stage information at each iteration $\nu$ into one new cut only, reducing the size of the master problem, especially in cases where $|\Omega|$ is large, at the expense of more iterations until the algorithm convergences.

 Finally, we remark that in both versions of the Benders algorithm, each subproblem \eqref{mod:suprolem} is essentially a different $\omega$-indexed parametrization of the same linear program. Therefore, no information exchange is required between the different instances of \eqref{mod:suprolem}, once the first-stage decisions are fixed using either the multi-cut version \eqref{mod:Benders-Master-MC} or the single-cut formulation \eqref{mod:Benders-Master-SC} of the master problem.


\vspace{-10pt}
\section{Data-driven adaptive Benders decomposition}
\label{sec:Improvements}
As previously discussed, the multi- and single-cut versions of the Benders decomposition algorithm outperform each other in terms of the number of iterations and the size of the master problem, respectively. In order to bridge this efficiency gap, we propose an adaptive data-driven version of the Benders algorithm that lies in-between the two prevailing approaches and aims at finding the optimal trade-off between the size of the master problem and the number of iterations.
This is a generic algorithm involving two distinct processes that can be applied independently to any problem that can be solved using Benders decomposition. The first process exploits the statistical information from the second-stage subproblems in order to make efficient use of the Benders cuts. The second process builds, on the other hand, an outer parallelization of the Benders algorithm which is based on first-stage information.  
In both processes, the $|\Omega|$ subproblems are treated independently and their optimal solutions are computed in parallel within each Benders iteration.

\vspace{-12pt}
\subsection{Improving Benders decomposition algorithm using data-driven techniques}

Our first data-driven process comprises two main functions: an intelligent cuts' aggregation based on data clustering techniques and the adaptive management of the existing and the new cuts in the course of the Benders algorithm. 
Aiming for a trade-off between the single- and multi-cut versions of the Benders algorithm, we group the scenarios into several clusters $c \in \mathcal{C}_k$ and we reformulate the master problem as:

\vspace{-12pt}  
\begin{subequations}
\begin{align}
  \underset{\Xi^{\text{DA}},\theta_{\omega} }{\text{Minimize}} & \quad C^{\text{DA}} + \sum_{\omega \in \Omega} \pi_{\omega} \theta_{\omega}
  \end{align}  
subject to
\begin{align}
& \text{constraints }  \eqref{eq:UC-FScns}, \; \eqref{eq:MC-cuts}, \\
& \theta_{\omega} \geq \theta^{\text{min}}, \quad \forall\omega, \\
& \!\!\! \sum_{\omega \in \Omega_c }  \pi_\omega\theta_\omega \ge \sum_{\omega \in \Omega_c }\pi_\omega \Theta_\omega^{(k)} \; : \; \mu_{ck}, \nonumber \\
& \quad \quad \quad \quad \quad \quad \forall k=\{1,...,\nu-1\} \setminus \mathcal{I}^\kappa, \forall c \in \mathcal{C}_k, \label{eq:Aggr-Cut} \\
& \!\!\! \sum_{\omega \in \Omega} \pi_{\omega} \theta_{\omega} \ge \sum_{\omega \in \Omega} \pi_{\omega} \Theta_{\omega}^{(k)}, \quad \forall k \in \mathcal{I}^{\kappa}, \label{eq:Aggr-Cons}
\end{align}
\end{subequations}
where $\Omega_c$ comprises all scenarios $\omega$ that belong to cluster $c$.

Using constraint \eqref{eq:Aggr-Cut}, we generate an aggregated cut for all scenarios $\omega \in \Omega_c$ and for each cluster $c \in \mathcal{C}_k$. This cut aggregation is applied to all cuts that were generated in each of the previous iterations $\{1,...,\nu-1\} \setminus \mathcal{I}^\kappa$, where $\mathcal{I}^\kappa$ is a dynamic set that contains the iteration counters of the inactive cuts for which we apply the cut consolidation technique from \cite{WOLF}. The cut consolidation technique performs an intelligent handling of the inactive cuts once they become obsolete. Indeed, only a limited number of cuts are binding at every iteration and usually many cuts remain inactive after certain iterations \cite{no_cut_removal_tech}.
In order to take advantage of this observation, 
we use constraint \eqref{eq:Aggr-Cons} to replace with a single cut all cuts that were created in iteration $k$ using the corresponding clusters $\mathcal{C}_k$ and remained inactive in $\kappa$ successive iterations, i.e. the respective Lagrange multiplier $\mu_{ck}$ of Benders cut \eqref{eq:Aggr-Cut} is equal to zero in $\kappa$ successive iterations. This strategy keeps the information from the inactive cuts in a compressed form and reduces the size of the master problem. It should be noted that the cut consolidation is applied only if all cuts generated at iteration $k$, using clusters $c \in \mathcal{C}_k$, are inactive for $\kappa$ successive iterations. Then iteration counter $k$ is appended to the set $\mathcal{I}^\kappa$ and the master problem constraints \eqref{eq:Aggr-Cut} and \eqref{eq:Aggr-Cons} are rebuilt according to the corresponding indices.

A relevant question that emerges after this reformulation of the master problem is how to define the set of clusters $\mathcal{C}$ such that the information contained in the solution of each $\omega$-indexed subproblem is passed as intact as possible to the master problem despite reducing the number of cuts. To tackle this question, we need to decide which clustering technique to use for computing the clusters and based on which attributes we will compare and evaluate  the similarity of the scenarios.

For every scenario, a number of attributes can be used to perform the comparison required to form the clusters. The most intuitive attribute is the production forecast $W^*_{jt\omega}$ of every wind farm $j$ and time $t$, which is the only parameter that varies in each subproblem for different scenarios $\omega \in \Omega$.
Considering that the wind power scenario set is a fixed parameter for the SUC problem, the cluster formation only needs to be performed once at the beginning of the Benders algorithm. 
However, this \textit{static} clustering approach is unable to account for any new information that may become available in the course of the Benders iterations, as the master problem traverses the second-stage value function.


In fact, at the end of each new Benders iteration $\nu$ we obtain updated values of the objective function \eqref{eq:obj-MC} and of the dual variables corresponding to constraints \eqref{eq:FixCns1}-\eqref{eq:FixCns4}. Both these values reflect the impact of first-stage decisions on the recourse cost for each  uncertainty realization $\omega$
and are the building blocks of the optimality cuts that drive the solution of the problem, i.e. similar values of these duals and the objective function \eqref{eq:obj-MC} would yield identical Benders cuts. To extract additional value from these data, we can apply a \textit{dynamic} clustering approach, using the aforementioned quantities as clustering attributes and re-compute the set of clusters $\mathcal{C}_k$ at every Benders iteration $k$.
   Four dual variables are defined for every scenario at every iteration, corresponding to different first-stage decisions, i.e. day-ahead wind power dispatch, power flows and upward/downward reserve procurement capacities. 
   In order to use these dual variables as clustering attributes, we need to normalize their values such that we can perform a meaningful comparison. As an illustration, the normalized value $\overline{\lambda}^{+}_{gt\omega}$ of the dual variable $\lambda^{+}_{gt\omega}$ is calculated as:
\begin{equation}\label{eq:norm}
    \overline{\lambda}^{+}_{gt\omega}=\ddfrac{\lambda^{+}_{gt\omega}-\min_{\forall g, t, \omega}(\lambda^{+}_{gt\omega})}{\max_{\forall g, t, \omega}(\lambda^{+}_{gt\omega})-\min_{\forall g, t, \omega}(\lambda^{+}_{gt\omega})}, \quad \forall g, t, \omega.
\end{equation}
Similarly, we obtain the normalized values of the remaining dual variables $\lambda^{-}_{gt\omega}$, $\lambda^{\text{W}}_{jt\omega}$, and $\lambda^{\text{F}}_{\ell t\omega}$.

Aggregation of Benders cuts into classes characterized by similar properties is paramount for obtaining an accurate description of the recourse value function, while restricting the number of constraints in the master problem. In general, clustering techniques are based on the comparison of objects with respect to a specific metric.
In Section \ref{sec:Results}, we evaluate and compare the performance of three clustering methods, namely, the k-shape, k-means and hierarchical algorithms. The k-shape \cite{Paparrizos} focuses primarily on the shape and the variation pattern of the data instead of their specific values. The k-means algorithm \cite{kaufman2009finding}, which belongs to the family of partitional clustering techniques, groups the data around k-centroids, where these centroids are the mean values of the scenarios within each k-cluster. Finally, in the agglomerative hierarchical method \cite{hastie2009unsupervised}, each object is initially considered as a cluster that are iteratively fused in pairs with minimum distance until the desired number of clusters is formed.

In order to balance efficiently the size of the master problem with the number of iterations, we develop an adaptive approach that varies the number of clusters until convergence is achieved.
The main idea is to introduce a measure for`convergence speed' and modify the number of clusters, considering that adding more clusters (cuts) improves the rate of convergence at consecutive iterations but increases the size of the master problem.
To quantify convergence speed, we utilize the property of the lower bound $\underline{\xi}$ being a non-decreasing function of the Benders iteration $\nu$. This can be intuitively verified considering that the master problem in iteration $(\nu+1)$ is a more constrained version of the $\nu$-th instance, due to the addition of the new cuts. Based on the above, we calculate the difference between two successive values of the lower bound as $\Delta=\underline{\xi}^{(\nu+1)}-\underline{\xi}^{(\nu)}$, which is then compared with two user-defined limits denoted as $\Delta^{\uparrow}$ and $\Delta^{\downarrow}$. If $\Delta < \Delta^{\uparrow}$ the number of clusters is increased by $\varrho$ clusters and vice versa if $\Delta > \Delta^{\downarrow}$. 
To make the algorithm less sensitive to small variations of $\Delta$, these thresholds are selected as $\Delta^{\uparrow} < \Delta^{\downarrow}$, defining a dead-band in which the number of clusters remains unchanged.

The complete algorithm (Algorithm 1) of the data-driven Benders decomposition methodology described above is provided in the Appendix of the paper.


\vspace{-8pt}
\subsection{Outer parallelization of the Benders algorithm}

Inspired by the work in \cite{Blanco}, we propose an outer parallelization process that comprises two steps. At the first step, we split the scenario set $\Omega$ into smaller subsets $s \in \mathcal{S}$ and we solve the SUC for all subsets.
Then, we compare the first-stage binary variables $u_{gt}$ at the optimal solution of each subset and we resolve the SUC problem using the full scenario set $\Omega$, fixing though the commitment variables that are identical in all subsets to the  solution obtained in the previous step. The complete algorithm (Algorithm 2) for the outer parallelization of the Benders algorithm is presented in the Appendix of the paper.

To form the set of subsets $\mathcal{S}$, we first apply the k-medoid clustering method to the original scenario set $\Omega$ in order to obtain $e \in \mathcal{E}$ clusters and the corresponding most significant scenario (medoid), denoted as $m_e$. Each subset $s$ is then constructed as a tuple $s$ = \{$\omega \in \Omega_e$\} $\bigcup$ $\{m_{e'} \; \forall e' \in \mathcal{E} \setminus e  \}$ that contains the scenarios of  cluster $e$ and the most significant scenarios of all the other clusters $e' \in \mathcal{E} \setminus e$. This subset formation approach is illustrated in Fig. A1 of the Appendix.
It is worth mentioning that the k-medoid  clustering technique is chosen here to ensure that the center of each cluster (medoid) corresponds to a scenario contained in $\Omega$ and not to an artificial trajectory as in the k-shape or k-means algorithms. As a result, each scenario subset $s$ consists of scenarios that have similar properties, i.e. belong to the same cluster $e$, as well as the medoids $m_{e'}$ from the rest of the clusters $e' \in \mathcal{E} \setminus e$. These scenarios $m_{e'}$ enrich the information contained in each subset, since they encode the general characteristics of uncertainty in their clusters without carrying along the particularities of all scenarios, such that the optimal solutions for different subsets can still be different.
To ensure notational clarity, we underline that the clusters $e \in \mathcal{E}$ used in the outer parallelization process differ from clusters $c \in \mathcal{C}_k$ employed in the cut aggregation process of the previous subsection. In particular, the set of clusters $\mathcal{E}$ is generated once, at the beginning of the outer parallelization process, based on the wind power scenario set $\Omega$, whereas the set of clusters $\mathcal{C}_k$ is re-computed in each iteration $k$ of the Benders algorithm.

The logic behind the process described above follows from the nature of the problem at hand. In particular, the primary purpose of the unit commitment problem is to find the optimal day-ahead schedule of the conventional units, in order to ensure that the system will have enough flexibility to cope with the variability and uncertainty of renewables during real-time operation. From a technical perspective\footnote{We acknowledge that the day-ahead schedule may have also economic implications, if the unit commitment model is used as market-clearing mechanism and the resulting production schedule is used in economic settlements.}, the rest of the decision variables provide rather an indicative/advisory dispatch schedule that is anyhow subject to changes, depending on the actual operating conditions. 
In this regard and considering that by construction the subsets describe adequately the diversity of plausible operating conditions, we  postulate that the common commitment variables among all subsets characterize accurately the true optimal solution. 

The benefits from applying this outer parallelization process is twofold. First, the SUC problem instances for each subset are smaller than the initial model with the complete scenario set. In addition, each SUC per subset can be solved independently and in parallel, exploiting at the same time all the data-driven techniques that we proposed in the previous subsection to speed up the Benders algorithm. An additional advantage is that the SUC problem run in the second step, despite including all scenarios in $\Omega$, is less computationally  intense, since a significant portion of the binary variables is fixed. Finally, this procedure provides the possibility to discard time-consuming subsets as we show in the following section. 


\section{Case Studies}
\label{sec:Results}

In this section, we evaluate the proposed data-driven adaptive version of the Benders algorithm and the outer parallelization process described in Section \ref{sec:Improvements}. We used the IEEE RTS-24 from \cite{ordoudisRTS-24} to evaluate the performance of the various clustering parameters, i.e. attributes and clustering techniques, and the IEEE 3-Area RTS-96 provided in \cite{RTS-96} to assess the efficiency  of the proposed improvements compared to the standard versions of the Benders decomposition algorithm. 
The IEEE RTS-24 system comprises 24 nodes connected by 34 transmission lines, 12 generators and 6 wind farms. The IEEE 3-Area RTS-96 consists of 72 nodes, 107 lines, 96 generators and 15 wind farms. 
We model wind power uncertainty using a set of scenarios that respect the spatio-temporal correlation of forecast errors over 15 different wind farm locations in Western Denmark. These scenarios are provided online at \cite{WindData}.

The tests on the IEEE RTS-24 were performed on a desktop computer, with an Intel(R) Xeon Gold(TM) 6154 CPU with 2 processors clocking at 2.99GHz and 479 GB of RAM. 
For the numerical experiments carried out on the IEEE 3-Area RTS-96, we used the SGE Arton Grid of D-ITET ETH Z\"urich \cite{ArtonGrid} with 11 computing servers of 16 cores Intel(R) Xeon(TM) at 2.9 GHz.
All test cases were implemented and solved in Python with Gurobi \cite{gurobi} as the MILP solver. The parallelization of the subproblems was implemented using the Joblib Python library \cite{joblib} with a multi-processing scheme. 
The input data and the corresponding Python codes are provided online in \cite{onlineapp}.

\subsection{Assessment of clustering parameters}

We use the IEEE RTS-24 system to evaluate and compare different attributes and clustering techniques for our application. The computations are performed using 50 equiprobable scenarios over a scheduling horizon of 10 hours, with stopping criterion $\varepsilon=10^{-6}$. 
In order to get an unbiased evaluation of the impact of different clustering parameters on computational time, the following tests are performed using the standard version of the Benders algorithm, without the improvements presented in Section \ref{sec:Improvements}.

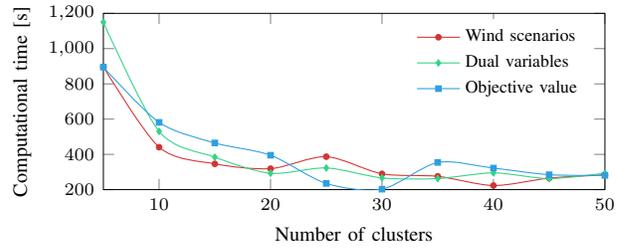
\begin{figure}[t]
  \begin{tikzpicture}[thick,scale=0.95]
\pgfplotsset{xmin=5, 
xmax=50,
ymin=200,
ymax=1200, 
try min ticks=5}
\begin{axis}[xlabel near ticks, ylabel near ticks, ylabel={Computational time [s]}, xlabel={Number of clusters}, label style={font=\footnotesize},
 tick label style={font=\scriptsize}, legend pos = north east, legend style={draw=none, name = Attr_leg_pos, font=\scriptsize, legend cell align={left}}, width=8.6cm, height=4.05cm]
    \addplot [color_ConvD!120, line width=0.15mm, mark=*, mark size=1.1, smooth, name path=A] 
    table[row sep=crcr]{
5 897.0126640796661\\
10 440.85258054733276\\
15 346.30490136146545\\
20 318.90883803367615\\
25 386.88130831718445\\
30 289.58475399017334\\
35 275.3513255119324\\
40 223.26364398002625\\
45 265.3455698490143\\
50 283.4931695461273\\
};
    \addlegendentry{Wind scenarios}
    \addplot [color_bilevel!120, line width=0.15mm, mark=diamond*, mark size=1.25, smooth] 
    table[row sep=crcr]{
5 1150.0750920772552\\
10 530.544605255127\\
15 384.0923891067505\\
20 292.7272605895996\\
25 323.0126304626465\\
30 266.5656051635742\\
35 262.87884855270386\\
40 295.00825572013855\\
45 262.6229989528656\\
50 293.2750380039215\\
};
    \addlegendentry{Dual variables}
    \addplot [color_StochD!120, line width=0.15mm, mark=square*, mark size=1.1, smooth] 
  table[row sep=crcr]{
5 894.8834369182587\\
10 581.5621018409729\\
15 464.7429177761078\\
20 395.10654950141907\\
25 235.01779079437256\\
30 202.5268771648407\\
35 354.0899660587311\\
40 322.1258280277252\\
45 284.56356596946716\\
50 280.4793257713318\\
};
    \addlegendentry{Objective value}
\end{axis}
\end{tikzpicture}
  \vspace{-8pt}
  \caption{Comparison of different clustering attributes for the SUC problem applied to the IEEE RTS-24 system, using the hierarchical clustering method.}
  \label{fig:compa_data_clust}
 \end{figure}

Figure \ref{fig:compa_data_clust} shows the computational time as a function of the number of clusters generated using the hierarchical method and considering as clustering attributes the wind power scenarios, the dual variables of the `fixing' constraints \eqref{eq:FixCns1} - \eqref{eq:FixCns4} and the subproblems' objective value. These different attributes exhibit similar performance. As the static clustering approach based on wind power scenarios is not compatible with an adaptive approach that varies the number of clusters between successive Benders iterations, we choose to proceed with the dynamic clustering based on dual variables, which according to Fig. \ref{fig:compa_data_clust} exhibits less performance variations compared to clustering based on the objective value. 
This follows from the fact that the latter clustering approach is based on a single-valued attribute per iteration, i.e. the objective function value, as opposed to clustering based on multiple dual variables that can encapsulate more information. 
 
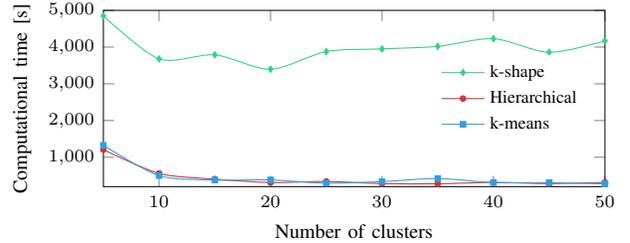
\begin{figure}[t]
\vspace{-12pt}
  \begin{tikzpicture}[thick,scale=0.95]
\pgfplotsset{xmin=5, 
xmax=50,
ymin=200,
ymax=5000, 
try min ticks=5}
\begin{axis}[xlabel near ticks, ylabel near ticks, ylabel={Computational time [s]}, xlabel={Number of clusters}, label style={font=\footnotesize},
legend style={at={(0.97,0.5)},anchor=east},
 tick label style={font=\scriptsize},  legend style={draw=none, name = Clust_leg_pos, font=\scriptsize, legend cell align={left}}, width=8.6cm, height=4.05cm]

    \addplot [color_bilevel!120, line width=0.15mm, mark=diamond*, mark size=1.25, smooth] 
  table[row sep=crcr]{
5 4849.039796829224\\
10 3678.4835221767426\\
15 3789.721709251404\\
20 3396.1364781856537\\
25 3878.535673379898\\
30 3950.7034063339233\\
35 4019.1150245666504\\
40 4229.398386001587\\
45 3857.0544719696045\\
50 4166.820436000824\\
};
\addlegendentry{k-shape}

    \addplot [color_ConvD!120, line width=0.15mm, mark=*, mark size=1.1, smooth, name path=A] 
  table[row sep=crcr]{
5 1208.223652601242\\
10 559.4292914867401\\
15 403.6804299354553\\
20 310.610915184021\\
25 340.55332922935486\\
30 281.37316727638245\\
35 277.9810469150543\\
40 314.86486053466797\\
45 282.5055196285248\\
50 313.68526458740234\\
};
\addlegendentry{Hierarchical}

\addplot [color_StochD!120, line width=0.15mm, mark=square*, mark size=1.1, smooth] 
  table[row sep=crcr]{
5 1325.0305325984955\\
10 499.11958932876587\\
15 379.4703297615051\\
20 382.03096199035645\\
25 292.60368824005127\\
30 336.8268995285034\\
35 419.55494689941406\\
40 309.55045533180237\\
45 307.60270524024963\\
50 274.51327419281006\\
};
\addlegendentry{k-means}

\end{axis}
\end{tikzpicture}
  \vspace{-8pt}
  \caption{Comparison of clustering methods for the SUC problem applied to the IEEE RTS-24 system, using dual variables as clustering attributes. }
  \label{fig:compa_clust_method}
  \vspace{-8pt}
\end{figure}

Moving ahead with clustering based on duals, we compare the performance of hierarchical, k-means and k-shape clustering techniques. The two first methodologies were implemented using the Scikit-learn Python library \cite{scikit-learn}, whereas for the third one we used the Sieve Python platform \cite{kshape}. Figure \ref{fig:compa_clust_method} shows that the k-shape method performs poorly for our problem compared to the hierarchical and k-means algorithms that can achieve notable reduction of computational time with few clusters. 
Hereinafter, we use the hierarchical clustering method due to its simpler implementation, without any compromise on the efficiency for our data-driven algorithm.

\subsection{Evaluation of data-driven Benders decomposition}

In this section, we use the IEEE 3-Area RTS-96 system to appraise the performance of the proposed data-driven improvements on Benders decomposition. The numerical tests are conducted for 100 scenarios and 24 hours. The MIP gap and the Benders convergence tolerance $\varepsilon$ are set to $10^{-6}$. 
To facilitate the assessment, we include in the algorithm one improvement scheme at a time and we compare the results in terms of computational time and problem size, expressed as the number of rows (constraints) of the master problem.

\subsubsection{Benders decomposition - standard versions}

Aiming to establish a benchmark for comparing our data-driven improvements, we first report in Table \ref{tab:para_SP} the results of the standard Benders decomposition algorithms, i.e. single- and multi-cut version. Both versions converge to the same objective function value, though the single-cut version is almost nine times slower than its multi-cut counterpart, affirming that the dense aggregation of recourse information in a single cut diminishes computational efficiency. On the other hand, the increased size of the master problem, despite being moderate in this particular case, is an indication of potential memory issues if the number of cuts grows significantly.

\begin{table}[t]
\vspace{-16pt}
\centering
\caption{Performance of single- and multi-cut Benders algorithms}
\begin{tabular}{@{}cccc@{}}
\toprule
  \begin{tabular}[c]{@{}c@{}}Benders\\ Algorithm\end{tabular}  & \begin{tabular}[c]{@{}c@{}}Expected\\ Cost [\$]\end{tabular}    & Time [s]  & \# of Rows \\ \midrule
 Single-cut & 747,007.74 & 62,074  & 23,292 \\ 
 Multi-cut  & 747,007.74 & 7,654 & 24,879  \\ \bottomrule
\end{tabular}
\vspace{-10pt}
\label{tab:para_SP}
\end{table}

\subsubsection{Adaptive cuts aggregation} The adaptive cuts aggregation strategy requires as input the number of clusters $\varrho$ that will be added if $\Delta$ lies outside the interval $[ \Delta^{\uparrow}, \Delta^{\downarrow} ]$. The initial number of clusters as well as the values of $\varrho$, $\Delta^{\uparrow}$, $\Delta^{\downarrow}$ are problem-specific parameters
and thus have to be tuned for the particular problem instance. In order to provide a more systematic approach for tailoring the thresholds $\Delta^{\uparrow}$ and $\Delta^{\downarrow}$, we compute them as percentages of the objective function value, i.e. $\mathcal{P}= \alpha \cdot [C^{\text{DA}} + C^{\text{BA}}]$, as follows:
\begin{align}
    \Delta^{\uparrow}=(1-\zeta)*\mathcal{P}, \quad \Delta^{\downarrow}=(1+\zeta)*\mathcal{P}, \label{eq:Delta}
\end{align}
where $\zeta$ sets a user-defined width of the dead-band around $\mathcal{P}$.

We set $\alpha=1 \%$ and we perform a grid search on the parameters $\zeta$ and $\varrho$. According to the results presented in Table \ref{tab:threshold_step}, we can see that the master problem grows larger as the cluster increment $\varrho$ increases, since more cuts are added in every iteration. Moreover,  the computational speed is optimized for $\varrho=5$, whereas there is no clear indication about the optimal value of $\zeta$. Based on these observations, in the following computations we set $\zeta = 75 \%$ and $\varrho= 5$, since this parameter setting provides a reasonable trade-off between the computational speed and the size of the master problem. The expected system cost is the same for all values of $\zeta$ and $\varrho$ and equal to the solution provided in Table \ref{tab:para_SP}.

\begin{table}[t]
\caption{Sensitivity analysis on cluster increment and threshold} 
\vspace{-4pt}
\centering
\begin{tabular}{@{}ccccc@{}}
\toprule
                                             &                                   & $\zeta = 25\%$ & $\zeta = 50\%$ & $\zeta = 75\%$ \\ \midrule
\multicolumn{1}{l|}{\multirow{2}{*}{$\varrho=1$}}  & \multicolumn{1}{l|}{Time {[}s{]}} & 4,959          & 5,024          & 6,400          \\
\multicolumn{1}{l|}{}                        & \multicolumn{1}{l|}{\# of Rows}   & 23,583         & 23,590         & 23,674         \\ \midrule
\multicolumn{1}{l|}{\multirow{2}{*}{$\varrho=5$}}  & \multicolumn{1}{l|}{Time {[}s{]}} & 4,100          & 4,211          & 4,096          \\
\multicolumn{1}{l|}{}                        & \multicolumn{1}{l|}{\# of Rows}   & 23,892         & 23,841         & 23,811         \\ \midrule
\multicolumn{1}{l|}{\multirow{2}{*}{$\varrho=10$}} & \multicolumn{1}{l|}{Time {[}s{]}} & 4,248          & 4,731          & 4,093          \\
\multicolumn{1}{l|}{}                        & \multicolumn{1}{l|}{\# of Rows}   & 24,060         & 24,144         & 24,001         \\ \bottomrule
\end{tabular}
\label{tab:threshold_step}
\vspace{-4pt}
\end{table}

Table \ref{tab:init_clust} shows the impact of the initial number of clusters, which along with the increment step $\varrho$ affects the number of cuts that are added in every Benders iteration. It becomes apparent that both the computational time and the number of rows increase when the initial number of clusters increases. Since we aim to minimize these two values, we choose to fix the initial number of clusters equal to 1.

\begin{table}[t]
\vspace{-16pt}
\centering
\caption{Sensitivity analysis on the initial number of clusters}
\begin{tabular}{@{}ccccc@{}}
\toprule
 \begin{tabular}[c]{@{}c@{}}Initial \\ \# of clusters\end{tabular}  & 1  & 25  & 50 & 100 \\\midrule
 Time [s] & 4,096 & 4,244  & 4,284  & 6,962\\ 
 \# of Rows& 23,811 & 23,864 & 24,054 & 24,789\\ \bottomrule
\end{tabular}
\vspace{-4pt}
\label{tab:init_clust}
\end{table}

\begin{table}[t]
\vspace{-8pt}
\centering
\caption{Sensitivity analysis on the number of subsets}
\vspace{-4pt}
\begin{tabular}{@{}ccccc@{}}
\toprule
 \# of Subsets $s$  & 1  & 10  & 30 \\ \midrule
 Exp. Cost [\$] & 747,007.74 & 747,359.74  & 747,007.74\\ 
 $T_1$ [s] &  \text{-}  &   15,819  & 3,697\\
 $T_2$ [s] &  \text{-}  &   1,452   & 1,391\\
 $T$   [s]    & 4,096 & 17,271  & 5,088  \\ 
 Max. \# of rows& 23,811 & 23,746 & 23,767\\ \bottomrule
\end{tabular}
\label{tab:sub_set}
\vspace{-8pt}
\end{table}

\begin{table}[h!]
\centering
\caption{Sensitivity analysis on the percentage of completed subsets}
\vspace{-4pt}
\begin{tabular}{@{}ccccc@{}}
\toprule
$\gamma$ & 100\%   & 85\% & 75\% & 25\% \\ \midrule
 Exp. Cost [\$] & 747,007.74  &747,007.74   & 747,007.74 & 747,130.28\\
 $T_1$ [s] & 3,697     &   2,657   &   2,434   & 2,199\\
 $T_2$ [s] & 1,391     &   1,350   & 1,364 & 1,518\\
 Max. \# of rows &23,767  & 23,671 & 23,671  & 24,098 \\ \bottomrule
\end{tabular}
\label{tab:share}
\vspace{-16pt}
\end{table}

\subsubsection{Outer parallelization strategy}

To assess the performance of this strategy we have to consider the total computational time $T = T_1 + T_2$, with $T_1$ and $T_2$ referring respectively to the time required for the parallel solution of subproblems $\omega \in \Omega_s$ within each subset $s$ and to the solution of SUC with fixed commitment variables. Denoting as $\tau_s$ the computational time of subset $s$, we define $T_1$ as $T_1=\max_{\forall s}(\tau_s)$.

Table \ref{tab:sub_set} compares the results of the SUC with (number of subsets $|\mathcal{S}|$ greater than 1) and without (number of subsets $|\mathcal{S}|$ equal to 1) the application of the outer parallelization strategy. 
Setting the number of subsets to $|\mathcal{S}|=10$ over-compresses and dilutes the information such that both the computational time and the final objective value are significantly increased. 
On the contrary, using the outer approximation strategy with $|\mathcal{S}|=30$, we manage to obtain the true optimal solution reducing at the same time the size of the problem. Indeed, the full benefits of problem size reduction may unfold when the SUC is applied to real-life power systems, when memory requirements become an issue. 
For $|\mathcal{S}|=30$, the algorithm computes six different commitment schedules (Commitment I - VI) as shown in Fig. A2 of the Appendix, which compares the commitment status in the subset solution with the corresponding optimal value.Only the points for which the binary commitment variables $u_{gt}$ that are different among the commitments and therefore not fixed in the second step of the outer parallelization strategy are shown. It is worth noting that each commitment schedule I, IV, V and VI is determined by one subset only, schedule II is given by two subsets and the remaining 24 subsets all yield the commitment schedule IV as the solution. This indicates that the commitment schedule is driven primarily by the uncertainty characteristics of each subset and not by the number of scenarios per se.

Based on these observations and taking into account that each subset in the outer parallelization strategy is solved independently and in parallel, we consider the option to disregard the most computationally demanding subsets. Table \ref{tab:share} reports the optimal value of the objective function, the solution time and the problem size that correspond to different percentages $\gamma$ of completed subproblems, e.g. $\gamma = 25 \%$ means that only a quarter of the subsets is solved to optimality before the first step of the outer parallelization algorithm is terminated. 

As already mentioned some commitment schedules are obtained only by one subset. Therefore, by reducing $\gamma$, we bear the risk to disregard unique schedules if their computational time is high. In that case, some commitment variables will be fixed to non-optimal values in the second step of the process, as if they were common for all subsets in $\mathcal{S}$. 
Having, however, only six different commitment schedules resulting from the full subset $\mathcal{S}$, we are indeed able  to recover  the true optimal solution even for $\gamma = 75\%$, albeit with reduced computational time and problem size.
In practice, this setting does not discover commitment schedules I and II within the $75\%$ of the subsets that are solved first and therefore fixes some extra commitment variables compared to the $\gamma = 100\%$ case that keeps them free. Nevertheless, in the second run of the SUC for $\gamma = 100\%$, these commitment variables obtain eventually the same optimal value as in the $\gamma = 75\%$ case where they were treated as fixed parameters.
This behavior changes for $\gamma = 25\%$, since this small subset of solutions does not capture adequately the different operating conditions. Only schedules III and VI are disclosed and five commitment variables are fixed to non-optimal values, resulting in increased expected system cost. It should be noted that we explored the possibility to  warm start instead of fixing the binary variables $u_{gt}$ in order to mitigate the risk of being trapped to sub-optimal solutions. However, this approach almost doubles the computational time $T_2$ and diminishes the advantages of the outer parallelization strategy.
 
\subsubsection{Cut Consolidation}

We apply the cut consolidation technique on the second-step SUC problem of the outer parallelization, which involves the complete set of uncertainty scenarios $\Omega$. 
Table \ref{tab:cut_conso} reports the computational time for different values of the iteration threshold $\kappa$, showing that this technique can reduce moderately the computational time $T_2$ and has a more pronounced impact on the size of the problem that is now reduced almost to the size of the single-cut version.

\begin{table}[H]
\centering
\vspace{-8pt}
\caption{Sensitivity analysis on the iteration threshold }
\vspace{-4pt}
\begin{tabular}{@{}ccccc@{}}
\toprule
 $\kappa$ & 2  & 5  & 10 \\ \midrule
 Exp. Cost [\$] & 747,007.74  & 747,007.74  & 747,007.74\\ 
 $T_2$ [s]&1,357 & 1,336 & 1,402\\
 \# of rows& 23,636 & 23,541 & 23,707\\ \bottomrule
\end{tabular}
\vspace{-8pt}
\label{tab:cut_conso}
\end{table}

\subsubsection{Performance analysis}

Table \ref{tab:compa} summarizes the performance of the different improvements and compares them to the extensive form of the SUC problem. All problem instances reach the same optimal solution.
The single-cut version of Benders algorithm has the worst performance in terms of computational time, whereas the multi-cut version can speed up considerably the solution process at the expense of increased problem size. 
Applying the adaptive clustering strategy can alleviate this issue and further reduce the computational time.
 The outer parallelization algorithm leverages the capabilities of distributed optimization to halve the computational time compared to the multi-cut approach, whereas cut consolidation allows to shrink the problem and almost reach the size of the single-cut version.
Overall, our data-driven version of Benders decomposition, achieves a two-fold reduction of computational time and reduces the problem size by a factor of 50 compared to the extensive form of the SUC problem.  

\begin{table}[H]
\vspace{-8pt}
\caption{Comparison of the improvements}
\vspace{-4pt}
\centering
\begin{tabular}{@{}llll@{}}
\toprule
  & Exp. Cost [\$] & Time [s]  & \# of rows  \\ \midrule
 Single-cut &747,007.74 & 62,074 & 23,292\\
 Multi-cut &747,007.74 & 7,654 & 24,879\\
 Cut aggregation & 747,007.74 & 4,096 & 23,811\\
 Outer parallelization & 747,007.74& 3,798 & 23,671\\
 Cut consolidation & 747,007.74& 3,770& 23,541 \\\midrule
 Extensive form &747,007.74& 8,477 & 1,402,779  \\ 
 \bottomrule
\end{tabular}
\label{tab:compa}
\vspace{-8pt}
\end{table}


\section{Conclusions}
\label{sec:conclusion}

In this paper, we propose several data-driven improvements of Benders decomposition in order to improve its computational performance and reduce its memory requirements, and we apply our algorithms to the stochastic unit commitment problem.
Using state-of-the-art clustering techniques we developed an adaptive cuts aggregation strategy that halves computational time compared to the multi-cut Benders approach. In addition, we proposed a novel outer parallelization approach that is implemented in a distributed fashion and is combined with a cut consolidation method. Our numerical results show that this data-driven scheme outperforms, in terms of computational time and problem size, both the standard Benders decomposition and the extensive form solved by a commercial optimization solver.

Future research may focus on the development of more systematic processes for tuning the user-defined parameters of the various data-driven solution schemes as well as on the extension of our algorithm to multi-stage formulations.

\bibliographystyle{IEEEtran}
\bibliography{bibliography}

\appendix

\setcounter{figure}{0}
\renewcommand{\thefigure}{A\arabic{figure}}

This document serves as an Appendix for the paper ``Data-driven Adaptive Benders Decomposition for the Stochastic Unit Commitment Problem".

\vspace{8pt}

\subsection{Algorithm of the data-driven Benders decomposition}

\begin{algorithm}
\caption{ Data-driven Benders decomposition}\label{alg:DataBenders}
\begin{algorithmic}[1]
\Require $\varrho$, $\Delta^{\uparrow}$, $\Delta^{\downarrow}$
\State $\nu := 0$, 
\Repeat 
\State \textbf{Solve} Master problem
\State Return optimal solution $\bm{x}^{(\nu)}$
\State Compute lower bound $\underline{\xi}^{(\nu)}$
\State \textbf{Set} $\bm{x}^{\text{sub}}_{\omega} := \bm{x}^{(\nu)}$ 
\State \textbf{Solve} in parallel subproblem (3a)-(3f), $\forall \omega \in \Omega$
\State Compute upper bound $\overline{\xi}^{(\nu)}$
\State Compute $\Delta$ and define number of clusters $c \in \mathcal{C}$
\State Generate clusters $c \in \mathcal{C}$
\State Add aggregated Bender cuts in the Master problem
\State Apply cut consolidation
\State $\nu := \nu+1$
\Until{$ |\overline{\xi}^{(\nu)}-\underline{\xi}^{(\nu)} |\leq \epsilon$}
\end{algorithmic}
\end{algorithm}

In Algorithm 1, vector $\bm{x}^{\nu}$ contains the optimal values of the primal variables of the master problem at iteration $\nu$, which are passed to the Benders subproblems according to constraints (3c)-(3f), i.e. $\bm{x}^{\nu} = \{r^{+}_{gt}, r^{-}_{gt}, \forall g,t; w_{jt}, \forall j,t ; \hat{f}_{\ell t}, \forall \ell,t\}$. 


\subsection{Algorithm of the outer parallelization of the Benders algorithm}

\begin{algorithm}
\caption{ Outer parallelization of the Benders algorithm}\label{alg:OuterPar}
\begin{algorithmic}[1]
\State Generate clusters $e \in \mathcal{E}$
\State Generate subsets $s \in \mathcal{S}$
\For{\textbf{each} $s \in \mathcal{S}$} 
\State \begin{varwidth}[t]{\linewidth}\textbf{Solve} SUC model (1) using data-driven \par 
                        \hskip\algorithmicindent Benders decomposition $\forall \omega \in s$ \end{varwidth}
\vspace{0.01pt}
\hspace{18.0em}\smash{$\left.\rule{0pt}{1.8\baselineskip}\right\}\ \mbox{in parallel}$}
\State Return optimal solution $\bm{u}_s$
\EndFor
\State Fix common binary commitment variables in $\bm{u}_s, \; \forall s \in \mathcal{S}$
\State \begin{varwidth}[t]{\linewidth}\textbf{Solve} SUC model (1) using data-driven \par 
                        \hskip\algorithmicindent Benders decomposition $\forall \omega \in \Omega$ \end{varwidth}
\end{algorithmic}
\end{algorithm}

In Algorithm 2, vector $\bm{u}_{s}$ contains the optimal values of the commitment status binary variables $u_{gt}$ for all subsets $s \in \mathcal{S}$. 


\newpage


\subsection{Supplemental figures}


  \begin{figure}[h]
  \centering
  \scalebox{0.8}{\usetikzlibrary{decorations.pathreplacing}
\usetikzlibrary{shapes.multipart}
\begin{tikzpicture}[every text node part/.style={align=center}]
\node[text width=3cm] (a) at (0,1) {Scenario \\ Set $\Omega$};
\node[circle,draw=gray!30,fill=gray!30,minimum size=10pt]  at (0,0) {\scriptsize $\omega_1$};
\node[circle,draw=gray!30,fill=gray!30,minimum size=14pt] at (0,-1) {\scriptsize $\omega_2$};
\node[circle,draw=gray!30,fill=gray!30,minimum size=14pt] at (0,-2) {\scriptsize $\omega_3$};
\node[circle,draw=gray!30,fill=gray!30,minimum size=14pt] at (0,-3) {\scriptsize $\omega_4$};
\node[circle,draw=gray!30,fill=gray!30,minimum size=14pt] at (0,-4) {\scriptsize $\omega_5$};
\node[circle,draw=gray!30,fill=gray!30,minimum size=14pt] at (0,-5) {\scriptsize $\omega_6$};

\draw[dotted] (1,1.0) -- (1,-5.5);

\node[text width=3cm] (a) at (2,1) {Cluster \\ Set $\mathcal{E}$};

\node[circle,draw=color_ConvD,fill=color_ConvD,minimum size=14pt] at (2,0) {\scriptsize $\omega_1$};
\node[circle,draw=gray!30,fill=gray!30,minimum size=14pt] at (2,-1) {\scriptsize $\omega_2$};
\node[circle,draw=color_bilevel,fill=color_bilevel,minimum size=14pt] at (2,-2) {\scriptsize $\omega_3$};
\node[circle,draw=gray!30,fill=gray!30,minimum size=14pt] at (2,-3) {\scriptsize $\omega_4$};
\node[circle,draw=gray!30,fill=gray!30,minimum size=14pt] at (2,-4) {\scriptsize $\omega_5$};
\node[circle,draw=color_StochD,fill=color_StochD,minimum size=14pt] at (2,-5) {\scriptsize $\omega_6$};

\draw [decorate,decoration={brace,amplitude=10pt}]
(1.7,-1.4) -- (1.7,0.4)  ;

\draw [decorate,decoration={brace,amplitude=10pt}]
 (2.3,0.4) -- (2.3,-1.4) node [black,midway,xshift=0.7cm] 
{ $e_1$};

\draw [decorate,decoration={brace,amplitude=10pt}]
(1.7,-3.4) -- (1.7,-1.6)  ;

\draw [decorate,decoration={brace,amplitude=10pt}]
(2.3,-1.6) -- (2.3,-3.4) 
node [black,midway,xshift=0.7cm] 
{ $e_2$};

\draw [decorate,decoration={brace,amplitude=10pt}]
(1.7,-5.4) -- (1.7,-3.6)  ;

\draw [decorate,decoration={brace,amplitude=10pt}]
(2.3,-3.6) -- (2.3,-5.4) 
node [black,midway,xshift=0.7cm] 
{ $e_3$};

\draw[dotted] (3.5,1.0) -- (3.5,-5.5);

\node[text width=3cm]  at (6,1.1) {Subsets $\mathcal{S}$};


\node[circle,draw=color_ConvD,fill=color_ConvD,minimum size=14pt]  at (4.5,-0.5) {\scriptsize $\omega_1$};
\node[circle,draw=gray!30,fill=gray!30,minimum size=14pt]  at (5.5,-0.5) {\scriptsize $\omega_2$};
\node[circle,draw=color_bilevel,fill=color_bilevel,minimum size=14pt]  at (7,-0.5) {\scriptsize $\omega_3$};
\node[circle,draw=color_StochD,fill=color_StochD,minimum size=14pt] at (8,-0.5) {\scriptsize $\omega_6$};

\node[text width=3cm]  at (5,0.4) {\scriptsize Cluster $e_1$};
\draw[dotted] (6.25,0.0) -- (6.25,-1.0);
\node[text width=3cm]  at (7.5,0.4) {\scriptsize \shortstack{Most significant \\ scenarios of $e_2$, $e_3$} };

\draw [decorate,decoration={brace,amplitude=6pt}]
(3.9,-1.0) -- (3.9,0.0)  ;

\draw [decorate,decoration={brace,amplitude=6pt}]
 (8.6,0.0) -- (8.6,-1.0) node [black,midway,xshift=0.7cm] 
{ $s_1$};


\node[circle,draw=color_bilevel,fill=color_bilevel,minimum size=14pt]  at (4.5,-2.5) {\scriptsize $\omega_3$};
\node[circle,draw=gray!30,fill=gray!30,minimum size=14pt]   at (5.5,-2.5) {\scriptsize $\omega_4$};
\node[circle,draw=color_ConvD,fill=color_ConvD,minimum size=14pt]  at (7,-2.5) {\scriptsize $\omega_1$};
\node[circle,draw=color_StochD,fill=color_StochD,minimum size=14pt] at (8,-2.5) {\scriptsize $\omega_6$};

\node[text width=3cm]  at (5,-1.6) {\scriptsize Cluster $e_2$};
\draw[dotted] (6.25,-2.0) -- (6.25,-3.0);
\node[text width=3cm]  at (7.5,-1.6) {\scriptsize \shortstack{Most significant \\ scenarios of $e_1$, $e_3$} };

\draw [decorate,decoration={brace,amplitude=6pt}]
(3.9,-3.0) -- (3.9,-2.0)  ;

\draw [decorate,decoration={brace,amplitude=6pt}]
 (8.6,-2.0) -- (8.6,-3.0) node [black,midway,xshift=0.7cm] 
{ $s_2$};


\node[circle,draw=gray!30,fill=gray!30,minimum size=14pt]  at (4.5,-4.5) {\scriptsize $\omega_5$};
\node[circle,draw=color_StochD,fill=color_StochD,minimum size=14pt]   at (5.5,-4.5) {\scriptsize $\omega_6$};
\node[circle,draw=color_ConvD,fill=color_ConvD,minimum size=14pt]  at (7,-4.5) {\scriptsize $\omega_1$};
\node[circle,draw=color_bilevel,fill=color_bilevel,minimum size=14pt] at (8,-4.5) {\scriptsize $\omega_3$};

\node[text width=3cm]  at (5,-3.6) {\scriptsize Cluster $e_3$};
\draw[dotted] (6.25,-4.0) -- (6.25,-5.0);
\node[text width=3cm]  at (7.5,-3.6) {\scriptsize \shortstack{Most significant \\ scenarios of $e_1$, $e_2$} };

\draw [decorate,decoration={brace,amplitude=6pt}]
(3.9,-5.0) -- (3.9,-4.0)  ;

\draw [decorate,decoration={brace,amplitude=6pt}]
 (8.6,-4.0) -- (8.6,-5.0) node [black,midway,xshift=0.7cm] 
{ $s_3$};

\end{tikzpicture}}
  \caption{Formation of clusters $e\in \mathcal{E}$ and subsets $s \in \mathcal{S}$ used in the outer parallelization of the Benders algorithm. Scenarios $\omega \in \Omega$ in color denote the most significant scenarios of the corresponding cluster $e$, according to the k-medoid clustering method.}
  \label{fig:SUC_para_first}
 \end{figure}
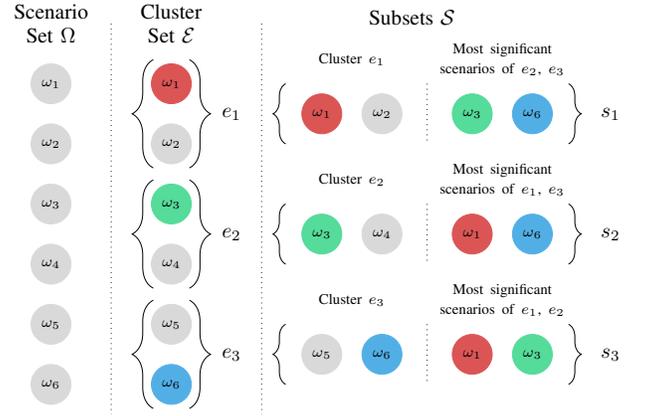

\pgfplotsset{every axis title/.append style={at={(0.5,0.85)}}}
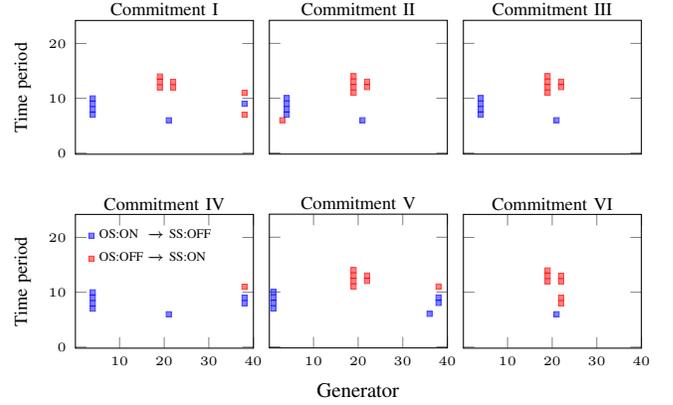
\begin{figure}[h]
\centering
\resizebox{8.8cm}{!}{%
\begin{tikzpicture}[thick,scale=1]
\pgfplotsset{
every non boxed y axis/.append style={y axis line style=-},
every non boxed x axis/.append style={x axis line style=-},
every tick label/.append style={font=\tiny},
compat=newest
}

\begin{axis}[
legend pos = north west,
legend style={fill=none, draw=none, font=\tiny, legend cell align={left}},
ylabel near ticks,
label style={font=\scriptsize},
ylabel style={align=center},
ylabel= {Time period},
title={\scriptsize Commitment I },
solid,
xmajorticks=false,
width=4cm,
height=3.4cm,
ymax = 24.2,
ymin  = -0.2,
xmin=0.1,
xmax=40,
scatter/classes={%
		a={mark=square*,blue,mark size=1pt,opacity=0.5},%
		b={mark=square*,red,mark size=1pt,opacity=0.5},%
		c={mark=o,draw=white,mark size=0pt}}]
]
\addplot [scatter,only marks,scatter src=explicit symbolic] table [x = G41, y = Time, meta = G41u]
{commit_schedule_diff_2.txt};
\addplot [scatter,only marks,scatter src=explicit symbolic] table [x = G191, y = Time, meta = G191u] {commit_schedule_diff_2.txt};
\addplot [scatter,only marks,scatter src=explicit symbolic] table [x = G211, y = Time, meta = G211u] {commit_schedule_diff_2.txt};
\addplot [scatter,only marks,scatter src=explicit symbolic] table [x = G221, y = Time, meta = G221u] {commit_schedule_diff_2.txt};
\addplot [scatter,only marks,scatter src=explicit symbolic] table [x = G381, y = Time, meta = G381u] {commit_schedule_diff_2.txt};
\end{axis}

\begin{axis}[
title={\scriptsize Commitment II},
solid,
xmajorticks=false,
yticklabels={},
width=4cm,
height=3.4cm,
xshift = 75,
ymax = 24.2,
ymin  = -0.2,
xmin=0.1,
xmax=40,
scatter/classes={%
		a={mark=square*,blue,mark size=1pt,opacity=0.5},%
		b={mark=square*,red,mark size=1pt,opacity=0.5},%
		c={mark=o,draw=white,mark size=0pt}}]
]
\addplot [scatter,only marks,scatter src=explicit symbolic] table [x = G32, y = Time, meta = G32u]
{commit_schedule_diff_2.txt};
\addplot [scatter,only marks,scatter src=explicit symbolic] table [x = G42, y = Time, meta = G42u] {commit_schedule_diff_2.txt};
\addplot [scatter,only marks,scatter src=explicit symbolic] table [x = G192, y = Time, meta = G192u] {commit_schedule_diff_2.txt};
\addplot [scatter,only marks,scatter src=explicit symbolic] table [x = G212, y = Time, meta = G212u] {commit_schedule_diff_2.txt};
\addplot [scatter,only marks,scatter src=explicit symbolic] table [x = G222, y = Time, meta = G222u] {commit_schedule_diff_2.txt};
\end{axis}

\begin{axis}[
title={\scriptsize Commitment III},
solid,
xmajorticks=false,
yticklabels={},
width=4cm,
height=3.4cm,
xshift = 150,
ymax = 24.2,
ymin  = -0.2,
xmin=0.1,
xmax=40,
scatter/classes={%
		a={mark=square*,blue,mark size=1pt,opacity=0.5},%
		b={mark=square*,red,mark size=1pt,opacity=0.5},%
		c={mark=o,draw=white,mark size=0pt}}]
]
\addplot [scatter,only marks,scatter src=explicit symbolic] table [x = G43, y = Time, meta = G43u]
{commit_schedule_diff_2.txt};
\addplot [scatter,only marks,scatter src=explicit symbolic] table [x = G193, y = Time, meta = G193u] {commit_schedule_diff_2.txt};
\addplot [scatter,only marks,scatter src=explicit symbolic] table [x = G213, y = Time, meta = G213u] {commit_schedule_diff_2.txt};
\addplot [scatter,only marks,scatter src=explicit symbolic] table [x = G223, y = Time, meta = G223u] {commit_schedule_diff_2.txt};
\end{axis}

\begin{axis}[
ylabel near ticks,
label style={font=\scriptsize},
ylabel style={align=center},
ylabel= {Time period},
legend pos = north west,
legend style={fill=none, draw=none, font=\tiny, legend cell align={left}},
title={\scriptsize Commitment IV},
solid,
width=4cm,
height=3.4cm,
yshift = -75,
ymax = 24.2,
ymin  = -0.2,
xmin=0.1,
xmax=40,
scatter/classes={%
		a={mark=square*,blue,mark size=1pt,opacity=0.5},%
		b={mark=square*,red,mark size=1pt,opacity=0.5},%
		c={mark=o,draw=white,mark size=0pt}}]
]
\addplot [scatter,only marks,scatter src=explicit symbolic] table [x = G44, y = Time, meta = G44u]
{commit_schedule_diff_2.txt};
\addlegendentry{OS:ON $\;\xrightarrow{}$ SS:OFF};
\addplot [scatter,only marks,scatter src=explicit symbolic] table [x = G214, y = Time, meta = G214u] {commit_schedule_diff_2.txt};
\addplot [scatter,only marks,scatter src=explicit symbolic] table [x = G384, y = Time, meta = G384u] {commit_schedule_diff_2.txt};
\addlegendentry{OS:OFF $\xrightarrow{}$ SS:ON};
\end{axis}

\begin{axis}[
xlabel={\footnotesize  Generator},
xlabel near ticks,
title={\scriptsize Commitment V},
solid,
yticklabels={},
width=4cm,
height=3.4cm,
xshift = 75,
yshift = -75,
ymax = 24.2,
ymin  = -0.2,
xmin=0.1,
xmax=40,
scatter/classes={%
		a={mark=square*,blue,mark size=1pt,opacity=0.5},%
		b={mark=square*,red,mark size=1pt,opacity=0.5},%
		c={mark=o,draw=white,mark size=0pt}}]
]
\addplot [scatter,only marks,scatter src=explicit symbolic] table [x = G15, y = Time, meta = G15u]
{commit_schedule_diff_2.txt};
\addplot [scatter,only marks,scatter src=explicit symbolic] table [x = G195, y = Time, meta = G195u] {commit_schedule_diff_2.txt};
\addplot [scatter,only marks,scatter src=explicit symbolic] table [x = G225, y = Time, meta = G225u] {commit_schedule_diff_2.txt};
\addplot [scatter,only marks,scatter src=explicit symbolic] table [x = G365, y = Time, meta = G365u] {commit_schedule_diff_2.txt};
\addplot [scatter,only marks,scatter src=explicit symbolic] table [x = G385, y = Time, meta = G385u] {commit_schedule_diff_2.txt};
\end{axis}

\begin{axis}[
title={\scriptsize Commitment VI},
solid,
yticklabels={},
width=4cm,
height=3.4cm,
xshift = 150,
yshift = -75,
ymax = 24.2,
ymin  = -0.2,
xmin=0.1,
xmax=40,
scatter/classes={%
		a={mark=square*,blue,mark size=1pt,opacity=0.5},%
		b={mark=square*,red,mark size=1pt,opacity=0.5},%
		c={mark=o,draw=white,mark size=0pt}}]
]
\addplot [scatter,only marks,scatter src=explicit symbolic] table [x = G196, y = Time, meta = G196u]
{commit_schedule_diff_2.txt};
\addplot [scatter,only marks,scatter src=explicit symbolic] table [x = G216, y = Time, meta = G216u] {commit_schedule_diff_2.txt};
\addplot [scatter,only marks,scatter src=explicit symbolic] table [x = G226, y = Time, meta = G226u] {commit_schedule_diff_2.txt};
\end{axis}
 \end{tikzpicture}
}
\vspace{-20pt}
\caption{Comparison between the commitment of the sub-sets and the optimal SUC's values. Blue/red points denote that the unit is ON/OFF in the optimal schedule (OS) of the SUC and OFF/ON in the subset solution (SS).}
\label{fig:commit_sub_set}
\vspace{-12pt}
\end{figure}

\end{document}